\def\QED{{\hfill$\Box$}}
\documentstyle[11pt,A4]{article}
\begin{document}
\bibliographystyle{plain}
\title{An asymptotic expansion inspired by Ramanujan%
\thanks{Appeared as Technical Report CMA-MR02-93/SMS-10-93, CMA, ANU,
February~1993.}\ %
\thanks{Copyright \copyright\ 1993--2010, R.~P.~Brent.
	\hspace*{\fill} rpb139tr typeset using \LaTeX}}
\author{Richard P.\ Brent\\
	Computer Sciences Laboratory\\
	Australian National University\\
	Canberra, ACT~0200}
\date{~}
\maketitle
\thispagestyle{empty}			%
\vspace{-1cm}
\begin{abstract}
Corollary 2, Entry 9, Chapter 4 of Ramanujan's first notebook
claims that
\[ \sum_{k=1}^\infty {(-1)^{k-1} \over nk} \left({x^k \over k!}\right)^n
	\;\sim\; \ln x + \gamma \]
as $x \to \infty$. This is known to be correct for the case $n = 1$,
but incorrect for $n \ge 3$. We show that the result is correct
for $n = 2$. We also consider the order of the
error term, and discuss a different, correct generalisation
of the case $n = 1$.

\smallskip

\smallskip


\end{abstract}

\section{Introduction}
\label{Sec:Intro}

Much of Ramanujan's work was not published during his lifetime,
but was summarised in his Notebooks. These were printed in facsimile
in 1957~\cite{N_fac}, and edited editions have been published
by Berndt~\cite{N1,N2,N3}.

Many of Ramanujan's results were obtained in a formal manner, and he did
not state sufficient conditions for their validity.  For example,
Example 3, Entry 9, Chapter 4 (page 97 of~\cite{N1}) is
\begin{equation}
\sum_{k=1}^\infty {(-1)^{k-1}(\varphi(k)-\varphi(-k)) \over k} = \varphi'(0)
								\label{eq:E3}
\end{equation}
but clearly some conditions on the function $\varphi(z)$ are required.
Berndt has given sufficient conditions, but they do not always hold for
Ramanujan's applications of~(\ref{eq:E3}).
In fact, Ramanujan does not claim exact equality in~(\ref{eq:E3}),
but writes that the left side is ``nearly'' equal to the right side.
Thus, some sort of approximation or asymptotic equality is intended.

To illustrate the use of~(\ref{eq:E3}),
take $\varphi(z) = x^z/\Gamma(z+1)$,
where $x > 0$ is real.
Proceeding formally, we obtain
\begin{equation}
\sum_{k=1}^\infty {(-1)^{k-1} \over k} {x^k \over k!}
	= \ln x + \gamma \;\;, 			\label{eq:right}
\end{equation}
where $\gamma = -\Gamma'(1)$ is Euler's constant.

If equality in~(\ref{eq:right}) is interpreted as asymptotic equality
(usually denoted by ``$\;\sim\;$'') as $x \to \infty$,
then the result is correct.
In fact, a classical result (also given on page 167 of~\cite{N2}) is
\begin{equation}
\sum_{k=1}^\infty {(-1)^{k-1} \over k} {x^k \over k!}
	\; - \; \ln x - \gamma =
	\int_x^\infty {e^{-t} \over t}dt = O(e^{-x}/x) \;\;. \label{eq:n1}
\end{equation}

Ramanujan's Corollary 2, Entry 9, Chapter 4 (page 98 of~\cite{N1})
is that, for positive integer $n$,
\[
\sum_{k=1}^\infty {(-1)^{k-1} \over k} \left({x^k \over k!}\right)^n
\;\sim\; n \sum_{k=1}^\infty {(-1)^{k-1} \over k} {x^k \over k!} \;\;.
\]
In view of~(\ref{eq:n1}), this is equivalent to
\begin{equation}
\sum_{k=1}^\infty {(-1)^{k-1} \over nk} \left({x^k \over k!}\right)^n
\;\sim\; \ln x + \gamma				\label{eq:C2}
\end{equation}
It is plausible that Ramanujan derived this result
from~(\ref{eq:E3}) in the same formal manner that
we derived~(\ref{eq:right}) above,
but taking $\varphi(z) = (x^z/\Gamma(z+1))^n$.

Berndt~\cite{N1} shows that~(\ref{eq:C2}) is false for $n \ge 3$;
in fact, the function defined by the left side of~(\ref{eq:C2})
changes sign infinitely often, and grows exponentially large as
$x \to \infty$. However, Berndt leaves the case $n = 2$ open.

The aim of this paper is to show that~(\ref{eq:C2}) is true in the
case $n = 2$. Theorem~\ref{Thm:T1} in Section~\ref{Sec:Main}
gives an exact expression for
the error in~(\ref{eq:C2}) as an integral involving the Bessel
function $J_0(x)$, and Corollary~\ref{Cor:C1} deduces an asymptotic
expansion. The most significant term is $O(x^{-3/2})$ as $x \to \infty$.

Theorem~\ref{Thm:T1} is a special case\footnote{Also given in
formula 11.1.20, Chapter 11 of Abramowitz and Stegun~\cite{Abr64}
(the chapter was written by Luke).}
of a formula given on page 48 of Luke~\cite{Luk62}.
However, the connection with Ramanujan does not seem to have
been noticed before.	%

In Corollary 2, Entry 2, Chapter 3 of his first Notebook,
Ramanujan shows that the function on the left side of~(\ref{eq:right})
can be written as
\[
e^{-x}\sum_{k=0}^\infty H_k{x^k \over k!} \;\;,	\label{eq:Hid}
\]
where $H_k = \sum_{j=1}^k {1/j}$ is a harmonic number ($H_0 = 0$).
Thus
\begin{equation}
\sum_{k=0}^\infty H_k{x^k \over k!} \biggl/
\sum_{k=0}^\infty    {x^k \over k!}
\;\sim\; \ln x + \gamma \;\;.
						\label{eq:rightH}
\end{equation}

In Section~\ref{Sec:Gen} we indicate how Ramanujan might have
generalised~(\ref{eq:rightH}) in much the same way that he
attempted to generalise~(\ref{eq:right}).

\section{Notation and Preliminary Results}
\label{Sec:Prelim}

In this section we give some preliminary results on integrals
involving $J_0$.
These results may be found in the literature, but for completeness
we sketch their proofs.

Recall that
\begin{equation}
J_0(x) = \sum_{k=0}^\infty {(-1)^k (x/2)^{2k} \over k!k!}\;\;. \label{eq:J0def}
\end{equation}
$J_0(z)$ is an entire function, but we are only concerned with its behaviour
on the positive real axis.  For small positive $x$,
$J_0(x) = 1 + O(x^2)$. For large positive~$x$, Hankel's
asymptotic expansion~\cite{Olv74,Wat44}		%
gives $J_0(x) = O(x^{-1/2})$.
These observations are sufficient to show that the integral
occurring in Lemma~\ref{Lemma:L2} below is absolutely convergent.

In the proof of the following Lemma,
$(a)_k = a(a+1) \cdots (a+k-1) =
\Gamma(a+k)/\Gamma(a)$, and		%
\[
F(a,b;c;z) = \sum_{k=0}^\infty {(a)_k (b)_k \over (c)_k}{z^k \over k!}
\]
is a hypergeometric function.

\newtheorem{lemma1}{Lemma}
\begin{lemma1}
For $0 < \mu < \frac{3}{2}$,
\[
\int_0^\infty t^{\mu - 1}J_0(t)dt =
{2^{\mu - 1}\Gamma(\mu/2) \over \Gamma(1 - \mu/2)} \;\;.
\]
\label{Lemma:L1}
\end{lemma1}
\noindent{\em Proof. }
The integral\footnote{More precisely, the generalisation with
$J_0(t)$ replaced by $J_\nu(t)$: see Sec.\ 13.24 of Watson~\cite{Wat44}.}
is known as Weber's infinite integral~\cite{Web68}.
We sketch a proof in the case $\mu < \frac{1}{2}$,
which is all that is needed below.

Let $\alpha$ be a complex variable.
We first evaluate
\[
I(\alpha) = \int_0^\infty e^{-\alpha t}t^{\mu - 1}J_0(t)dt
\]
for $R(\alpha) > 1$. Integrating term by term,
the power series~(\ref{eq:J0def}) gives
\begin{equation}
I(\alpha) =
\alpha^{-\mu}\Gamma(\mu)
F\left({\frac{\mu}{2}},{\frac{\mu+1}{2}};1;-\alpha^{-2}\right)
\;\;.							\label{eq:ex5.7}
\end{equation}
By analytic continuation, the result~(\ref{eq:ex5.7}) holds\footnote{A
generalisation of~(\ref{eq:ex5.7}) is given in
Sec.\ 13.2 of Watson~\cite{Wat44}, and is attributed to
Hankel~\cite{Han75} and Gegenbauer~\cite{Geg76}.}
in the right half-plane $R(\alpha) > 0$.
Using a well-known hypergeometric function
identity\footnote{See equation (10.16), Chapter 5 of
Olver~\cite{Olv74}.}, %
we deduce that
\[
I(\alpha) =
{\Gamma(\mu) \over \Gamma(1 - {\frac{\mu}{2}})}\left(
{\Gamma({1 \over 2})F({\mu \over 2}, {\mu \over 2}; {1 \over 2}; -\alpha^2)
\over \Gamma({\mu + 1 \over 2})} +
{\alpha\Gamma(-{1 \over 2})F({\mu + 1 \over 2}, {\mu + 1 \over 2};
	{3 \over 2}; -\alpha^2)
\over \Gamma({\mu \over 2})}\right) \;\;.
\]
Since $J_0(t) = O(t^{-1/2})$ as $t \to \infty$,
our assumption that $\mu < \frac{1}{2}$ makes it is easy to justify
changing the order of integration
and taking the limit as $\alpha \to 0+$. Thus
\[
\int_0^\infty t^{\mu - 1}J_0(t)dt =
\lim_{\alpha \to 0+} I(\alpha) =
{\Gamma(\mu)\Gamma({1 \over 2}) \over
\Gamma(1 - {\mu \over 2})\Gamma({\mu + 1 \over 2})} \;\;.
\]
The Lemma now follows from the duplication formula
for the Gamma function.
  							\QED

\newtheorem{lemma2}[lemma1]{Lemma}
\begin{lemma2}
\begin{equation}
\int_0^\infty \left({e^{-t/2} - J_0(t) \over t}\right)dt = 0 \;\;.
\label{eq:int2}
\end{equation}
\label{Lemma:L2}
\end{lemma2}
\noindent{\em Proof. } A slightly more general result is given
in equation 6.622.1 of
Gradshteyn and Ryzhik~\cite{Gra65}, and attributed to
Nielsen~\cite{Nie06}.
We show that~(\ref{eq:int2}) follows easily from Lemma~1.

Let $\mu$ be a small positive parameter.
From Lemma~\ref{Lemma:L1},
\begin{equation}
\int_0^\infty (e^{-t/2} - J_0(t))t^{\mu - 1}dt =
2^\mu\Gamma(\mu) - {2^{\mu - 1}\Gamma(\mu/2) \over \Gamma(1 - \mu/2)} \;\;.
							\label{eq:diff}
\end{equation}
Since $\Gamma(1) = 1$, $\Gamma'(1) = -\gamma$, and
$\mu\Gamma(\mu) = \Gamma(1 + \mu)$, the right side of~(\ref{eq:diff}) is
\[
{2^\mu \over \mu}\left((1 - \gamma\mu) -
	\left({{1 - \gamma\mu/2} \over {1 + \gamma\mu/2}}\right) +
	O(\mu^2)\right) = O(\mu) \;\;.
\]
The result follows on letting $\mu \to 0+$.
							\QED

\section{Ramanujan's Corollary for $n = 2$}
\label{Sec:Main}

Our main result is the following Theorem,
which proves that~(\ref{eq:C2}) is valid for $n = 2$.

\newtheorem{thm1}{Theorem}
\begin{thm1}
Let
\[
e(x) = \sum_{k=1}^\infty {(-1)^{k-1} \over 2k}
	\left({x^k \over k!}\right)^2 - \ln x - \gamma \;\;.
\]
Then, for real positive $x$,
\[
e(x) = \int_{2x}^\infty {J_0(t) \over t}dt \;\;.
\]
\label{Thm:T1}
\end{thm1}
\noindent{\em Proof. } Proceeding as on page 99 of~\cite{N1},
using the fact\footnote{See page 103 of~\cite{N1}.
} that
\[
\gamma = \int_0^1 {1 - e^{-t} \over t} dt -
	\int_1^\infty {e^{-t} \over t} dt \;\;,
\]
we have
\[
e(x)  = \int_0^x {1 - J_0(2t) \over t} dt
	- \int_1^x {dt \over t}
	- \int_0^1 {1 - e^{-t} \over t} dt
	+ \int_1^\infty {e^{-t} \over t} dt
\]
\[      = \int_0^x {e^{-t} - J_0(2t) \over t}dt
	+ \int_x^\infty {e^{-t} \over t}dt \;\;.
\]
Now, from Lemma~\ref{Lemma:L2} with a change of variable,
\[
\int_0^x {e^{-t} - J_0(2t) \over t}dt
= \int_x^\infty {J_0(2t) - e^{-t} \over t}dt \;\;,
\]
so
\[
e(x) = \int_x^\infty {J_0(2t) - e^{-t} \over t}dt
	+ \int_x^\infty {e^{-t} \over t}dt
	= \int_x^\infty {J_0(2t) \over t} dt \;\;,
\]
and Theorem~\ref{Thm:T1} follows by a change of variable.
							\QED

\newtheorem{cor1}{Corollary}
\begin{cor1}
Let $e(x)$ be as in Theorem~\ref{Thm:T1}.
Then, for large positive $x$, $e(x)$ has an asymptotic expansion whose leading
terms are given by
\[
e(x) = {1 \over 2\pi^{1/2}x^{3/2}}
	\left(\cos \left(2x + {\frac{\pi}{4}}\right) +
	{13 \sin \left(2x + {\frac{\pi}{4}}\right) \over 16x} +
	O(x^{-2})\right) \;\;.
\]
\label{Cor:C1}
\end{cor1}
\noindent{\em Proof. } Using integration by parts and the fact that
$ xJ_0''(x) + J_0'(x) + xJ_0(x) = 0 $
(a special case of Bessel's differential equation),
it is easy to deduce from Theorem~\ref{Thm:T1} that
\[
e(x) = {J_0'(2x) \over 2x} + {J_0(2x) \over 2x^2}
	- 4\int_{2x}^\infty {J_0(t) \over t^3} dt \;\;.
\]
Continuing in the same way, we obtain an asymptotic expansion
\[
e(x) \;\sim\; {J_0'(2x) \over 2x} \sum_{k=0}^\infty (-1)^k {k!k! \over x^{2k}}
    \;+\; {J_0(2x) \over 2x^2}\sum_{k=0}^\infty (-1)^k {k!(k+1)! \over x^{2k}}
\;\;.\]
The Corollary follows from
Hankel's asymptotic expansions for $J_0(z)$
and $J_0'(z) = -J_1(z)$.
        						\QED

\section{A Different Generalisation}
\label{Sec:Gen}

An obvious generalisation of~(\ref{eq:rightH}) is
\begin{equation}
{\sum_{k=0}^\infty H_k\left({x^k \over k!}\right)^n \biggl/
{\sum_{k=0}^\infty \left({x^k \over k!}\right)^n}}
\;\sim\; \ln x + \gamma
						\label{eq:BM1}
\end{equation}
as $x \to \infty$.

It is easy to show that~(\ref{eq:BM1}) is valid for all positive
integer $n$. An essential difference
between~(\ref{eq:C2}) and~(\ref{eq:BM1})
is that there is a large amount of cancellation between terms
on the left side of~(\ref{eq:C2}),
but there is no cancellation in the
numerator and denominator on the left side of~(\ref{eq:BM1}).
The function $(x^k/k!)^n$ acts as a smoothing kernel with a peak at
$k \simeq x - 1/2$.	%
Since
\[
H_k = \ln k + \gamma + O(1/k) \;\;,
\]
the result~(\ref{eq:BM1}) is not surprising. What may be surprising is
the speed of convergence.
Brent and McMillan~\cite{RPB49} show that\footnote{We
note an error on page 310 of~\cite{RPB49}:
in the definition of $V_p(z)$,
``$z/k!$'' should be ``$z^k/k!$''.}	%
\begin{equation}
{\sum_{k=0}^\infty H_k\left({x^k \over k!}\right)^n \biggl/
{\sum_{k=0}^\infty \left({x^k \over k!}\right)^n}}
= \ln x + \gamma + O(e^{-c_nx})
						\label{eq:BM2}
\end{equation}
as $x \to \infty$, where
\[
c_n = \cases{	1,			&if $n = 1$;\cr
                2n\sin^2(\pi/n),	&if $n \ge 2$.\cr}
\]
In the case $n = 2$, (\ref{eq:BM2}) has error $O(e^{-4x})$.
Brent and McMillan
used this case with $x \simeq 17,400$ to compute
$\gamma$ to more than 30,000 decimal places.
From Corollary~\ref{Cor:C1}, the same value of $x$ in~(\ref{eq:C2})
would give less than 8-decimal place accuracy\footnote{More than
15,000 decimal places would have to be used in the computation to compensate
for cancellation of terms $\Omega_\pm(e^{2x}/x^2)$ in~(\ref{eq:C2})~!}.

The case $n = 3$ of~(\ref{eq:BM2}) is interesting because
$\max c_n = c_3 = 4.5$. However, no one seems to have used $n > 2$
in a serious computation of~$\gamma$.

Although Ramanujan~\cite{N2,R_works} gave many rapidly-convergent
series and other good approximations for~$\pi$,
he does not seem to have given series which are particularly
useful for approximating~$\gamma$, except for~(\ref{eq:n1})
and~(\ref{eq:rightH}) above.
In his paper~\cite{R_gamma} on series for~$\gamma$,
he gives several interesting series, of which the simplest\footnote{Due
to Glaisher: see~\cite{R_gamma}.} %
is
\[
\gamma = 1 - \sum_{k=1}^\infty {\zeta(2k+1) \over (k+1)(2k+1)} \;\;,
\]
but these series all involve the Riemann zeta function or related functions,
so they are not very convenient for computational purposes.

Our analysis has assumed that $n$ in
(\ref{eq:C2}) and~(\ref{eq:BM2}) is a positive integer.
It would be interesting to consider the behaviour of the functions
occurring in these equations for positive but non-integral
values of $n$, especially in the range $1 < n < 2$.

\end{document}